\documentclass[11pt]{article}
\usepackage{amsmath,amssymb,amsfonts,amsthm}
\usepackage{hyperref}
\usepackage[latin1]{inputenc}
\usepackage{epsfig}
\newtheorem{thm}{Theorem}

\newtheorem{prop}{Proposition}
\newtheorem{lem}{Lemma}
\newtheorem{assum}{Assumption}
\newtheorem{rem}{Remark}

\newtheorem{coro}{Corollary}

\newcommand{\dE}{\mathbb{E}}
\newcommand{\dN}{\mathbb{N}}
\newcommand{\dP}{\mathbb{P}}
\newcommand{\dR}{\mathbb{R}}

\newcommand{\ind}{\mathbf{1}}
\newenvironment{dem}[1][]%
   {\ \\ {\bf Proof #1. }}%
   {\hfill\mbox{\rule{2 true mm}{3 true mm}}}

   \title{Quantitative exponential bounds for the renewal theorem with
     spread-out distributions}

   \author{J.-B.~Bardet \footnote{Laboratoire de Math\'ematiques
       Rapha\"el Salem (LMRS), UMR 6085 CNRS, Universit\'e de
       Rouen. {Avenue de l'Universit\'e, BP 12, 76801 Saint Etienne du
         Rouvray, France} \
       E-mail:\texttt{jean-baptiste.bardet@univ-rouen.fr}} ,\,
     A.~Christen \footnote{Instituto de Estad\'istica, Pontificia
       Universidad Cat\'olica de Valpara\'iso. Errazuriz
       2734. Valparaíso, Chile.
       E-mail:\texttt{alejandra.christen@ucv.cl}} \, and \,
     J.~Fontbona\footnote{Department of Mathematical Engineering and
       Center for Mathematical Modeling, UMI(2807) UCHILE-CNRS,
       University of Chile, { Casilla 170-3, Correo 3, Santiago,
         Chile.}  E-mail:\texttt{fontbona@dim.uchile.cl}.  }}

\begin{document}
\maketitle

\begin{abstract} We establish exponential convergence estimates for
  the renewal theorem in terms of a uniform component of the
  inter-arrival distribution, of its Laplace transform which is
  assumed finite on a positive interval, and of the Laplace transform
  of some related random variable.  Although our bounds are not sharp,
  our approach provides tractable constructive estimates for the
  renewal theorem which are computable (theoretically and numerically,
  at least) for a general class of inter-arrival distributions. The
  proof uses a coupling, and relies on Lyapunov-Doeblin type arguments
  for some discrete time regenerative structure, which we associate with the
  renewal processes.

\end{abstract}

 

\section{Introduction and main statements}

We consider a classic {\it renewal processes} ${(T_n)}_{n\geq0}$
defined by $ T_n=T_0+\sum_{i=1}^n X_i\,, $ with $T_0$ a given
non-negative random variable, called {\it delay}, and
${(X_i)}_{i\geq1}$ an independent sequence of i.i.d.\ random
variables, equal in law to a given random variable $X>0$ with finite
mean. The random variables $(X_n)_{n\geq 1}$ and $(T_n)_{n\geq 0}$ are
respectively called {\it inter-arrivals times} and {\it renewal
  instants} (or epochs).  The {\it renewal measure} $U$, defined on
$(0,\infty)$ by $U(dt)= \dE \sum_{j=0}^{\infty}\delta_{T_j}(dt)$ (with
$\delta_x$ the Dirac mass at $x > 0$), is the central object of study
in renewal theory.  In one of its simplest forms, the {\it Renewal
  Theorem} states that, asymptotically as a time parameter $t>0$ goes
to infinity, the renewal measure of an interval $(t,t+h]$ is
proportional to $h>0$, if the distribution of $X$ is {\it
  non-arithmetic} (i.e.\ it is not supported on some real arithmetic
sequence). More precisely, one has
\begin{equation}
  \label{convergerenewal} 
  U((t,t+h])\to\frac{h}{\mu} \quad \mbox{ as }t\to \infty,
\end{equation}
where $0<\mu=\dE(X)< +\infty$. Originally established in the
non-arithmetic case in \cite{Black0} (and in \cite{ErdFellPoll} in the
arithmetic case), the Renewal Theorem and its several proofs have
motivated deep probabilistic ideas and developments. We refer to the
work \cite{Feller} for an analytic proof based on Choquet-Deny's Lemma
and to \cite{Lindvall} for the first probabilistic proof using
coupling, both ultimately relying on the Hewitt-Savage Theorem.
Self-contained probabilistic proofs were given in
\cite{Ney,ThorriComplete,LindvallRogers}.  See also
\cite{LindCoupling,kalashnikov94,alsmeyer2,alsmeyer,asmussen03} and
references therein for further background as well as for refinements
or extensions of the Renewal Theorem, and p.\ 480 in \cite{ThoriBook}
or the unpublished notes \cite{alsmeyerbook} for historical accounts.

It is well known that the tail of $X$ qualitatively determines the
asymptotic behavior of the renewal measure. For instance, if
$\dE(X^2)<\infty$ then the number of renewals in $(0,t]$ exhibits
Gaussian fluctuations as $t$ goes to infinity (see e.g.\ Prop.\ 6.3,
Ch.\ V in \cite{asmussen03}); if furthermore $X$ has some finite
exponential moment and is spread-out (see below for the context), the
error in \eqref{convergerenewal} is exponentially small (see e.g.
Thm.\ 2.10, Ch.\ VII in \cite{asmussen03}).  However, besides some
specific families of inter-arrival laws, the precise relation between
the tail of $X$ and the rate of convergence in the renewal theorem is
only partially understood. For instance, in \cite{BerenLund}, sharp
estimates were obtained in the arithmetic case, but only for monotone
hazard rates. Some conditions relating the renewal convergence rate to
the tail of $X$ in the arithmetic case are discussed in
\cite{Giacomin}. Estimates in the spread-out case have been given in
\cite{Ney}, \cite{asmussen03} but they depend on asymptotic bounds on
the renewal measure or equation.

The present note further explores the link between $X$ and the speed
of convergence in the renewal theorem by providing, for a wide class
of inter-arrival distributions, tractable bounds which are computable
(theoretically and numerically, at least) in terms of the law of $X$.

For the sake of concreteness, we will focus on inter-arrival
distributions which have some finite exponential moment and we will
furthermore assume they have a uniform component (which grants
non-arithmeticity and allows for a simpler analysis). More precisely,
introducing the notation
${\cal L}(\beta):= \dE\big(e^{\beta X}\big), \, \beta\in \dR$ for the
Laplace transform of $X$, we will make
\begin{assum}[exponential moment]
  \label{as:expo}
  The inter-arrival distribution admits some finite exponential
  moment:
  $\, \exists \, \alpha>0 \mbox{ such that }{\cal L}(\alpha
  )<+\infty.$
\end{assum}
We will also suppose that the law of $X$ satisfies
\begin{assum}[uniform component]
  \label{as:spread}
  There exist $c\geq L>0$ and $ \tilde{\eta}\in (0,1)$ such that
  \[
  \dP\big(X\in[t_1,t_2]\big)\geq\frac{\tilde{\eta}}{2L}
  \int_{t_1}^{t_2}\ind_{[c-L,c+L]}(r)\,dr \quad \mbox{ for all } 0\leq
  t_1\leq t_2 .
  \]
  In other words, $X$ has a uniform component on the interval
  $[c-L,c+L]$ with mass $\tilde{\eta}$.  
\end{assum}

In concrete examples, uniform components can usually be explicitly
identified. Recall also that Assumption \ref{as:spread} holds for some
convolution power of each spread-out distribution (i.e.\ one for which
some convolution power has an absolutely continuous component, see
Section VII.1 in \cite{asmussen03}).  Thus, our results also apply to
spread-out distributions by considering some finite sum of the
inter-arrival time lengths instead of a single one.

As in previous works, our approach will be based on a coupling
argument, that is, on constructing on some probability space two
copies of the renewal process with different initial delays, and
estimating the tail of some random time at which they ``coalesce''.
We briefly recall some general well known facts about such a
construction (see \cite{LindCoupling,ThoriBook,asmussen03} for more
background) and then state our results.

Write
$N_t:=\sup \{n\in \dN: T_n\leq
t\}=\sum_{j=0}^{\infty}\ind_{(0,t]}(T_j) $ for the total number of
renewals until time $t\geq 0$ and denote the {\it residual life} (or
forward recurrence time) process by
$$(B_t:=T_{N_t+1}-t)_{t\geq 0}\, ,$$ which is Markov.  Let
$(B_t')_{t\geq 0}$, $(N_t')_{t\geq 0}$ and $U'$ denote the
corresponding objects associated with a copy ${(T'_n)}_{n\geq 0}$ of
the renewal process with same inter-arrival law, defined on the same
probability space as ${(T_n)}_{n\geq 0}$.  An almost surely finite
random time $T^*$ such that a.s., $B_t=B_t'$ for all $t\geq T^*$, is
called a {\it coupling time for} $(B,B')$.
  
The recurrent process $(B_t)_{t\geq 0}$ has the stationary density
$\mu^{-1}\dP(X>t)dt$ and the renewal process with delay $T_0$
accordingly distributed is stationary (i.e.\ the corresponding renewal
measure is equal to $\mu^{-1}dt $).  The spread-out condition is
necessary and sufficient for the residual life process to converge in
total variation distance $\| \cdot \|_{TV}$ to its stationary
distribution (see e.g.\ Cor.1.5 Ch.VII in \cite{asmussen03}).  By the
coupling inequality (see \cite{LindCoupling}) one moreover has the
estimate
\begin{equation}
  \label{couplingineq} 
  \| law(B_t) - law(B_t') \|_{TV}  \leq \dP( T^*>t).
\end{equation}
Thus, finiteness of some exponential moment of $T^*$ immediately
grants exponential convergence to equilibrium at the same rate at
least, by Chernoff's inequality.  Moreover, since
$N_{T^* +s}-N_{T^*}= N'_{T^* +s}-N'_{T^*}$ a.s.\ for all $s\geq 0$,
for any given $t\geq 0$ one also gets the estimate
\begin{equation}
  \label{couplingrenew}
  \left| U(t+D)-  U'(t+D)\right|\leq \dE\left(\ind_{T^*>t}  
    \left(\sum_{j=0}^{\infty}\ind_{t+D}(T_j)  
      + \sum_{j=0}^{\infty}\ind_{t+D}(T'_j) \right)  
  \right) 
\end{equation}
for all Borel sets $D\subset\dR_+$.  Our  goal thus is to build
two copies ${(T_n)}_{n\geq 0}$ and ${(T'_n)}_{n\geq 0}$ with a
coupling time $T^*$ having an exponential tail that can be explicitly
controlled in terms of the law of $X$.
  
Let us introduce further notation required to state our results. In
the sequel we write
\[
\eta:=\tilde{\eta}/2\in \left(0,1/2\right)
\]
for the mass of the uniform component $[c,c+L]$ of $X$. We will also
denote by $\bar{{\cal L}}_a:\dR\mapsto \dR_+\cup \{\infty\}$ the
Laplace transform of the maximum of two independent copies of the
random variable $X$, both conditioned on being strictly larger than
$a>0$.  Last, given $x\geq 0$, we denote by $(B^x_t)$ the residual
life process when $T_0=x$ a.s.  The following is our main result:

\begin{thm}
  \label{main}
  Suppose Assumptions 1.\ and 2.\ hold. Given 
  $ \beta\in \dR_+, \delta \in [0,1)$ such that\\
  ${\cal L}( (1+\delta)\beta) <\infty $, set
  \[
  R=R(\delta,\beta):=\frac{1}{2\beta}\log \left[\frac{{\cal
        L}((1+\delta)\beta)}{ 1-{\cal L}(-(1-\delta) \beta) }\right].
  \]
  For each $ x\in \dR_+$, there exists a coupling 
  $(B_t^x,B_t^0)_{t\geq 0}$ with coupling time
  $T^*(x)$ such that
  \[
  \dP(T^*(x)>t) \leq \exp \left\{\theta \beta x \mathbf{1}_{x >R}
  \right\} \frac{\eta^{ \lceil R/L\rceil} e^{\theta \beta \left\{ R
        +\lfloor R/L\rfloor c\right\}} \bar{{\cal L}}_{c+L} \left(
      \theta \beta\right) }{1- e^{\theta \beta \left\{ R +\lfloor
        R/L\rfloor c\right\}} \bar{{\cal L}}_{c+L} \left( \theta
      \beta\right)(1-\eta^{ \lceil R/L\rceil} )}\, \exp\left( -\theta
    \delta \beta t \right)
  \]
 for every $ t\geq x$ and  all $\theta \in (0,1] $ for which
  $ e^{\theta \beta \left\{ R +\lfloor R/L\rfloor c\right\}}
  \bar{{\cal L}}_{c+L} \left( \theta \beta\right) \left(1-\eta^{
      \lceil R/L\rceil}\right)<1$.
\end{thm}

The above bound is involved, but can be better understood in terms of
the parameters $$a=(1+\delta)\beta \mbox{ and }b=(1-\delta)\beta.$$
Indeed, $R>0$ will correspond to the smallest value for which we can
grant that $R$-close renewals of two independent copies will occur
within some random time lapse $T_R>0$ of exponentially decaying
length. This value is controlled by both a positive and a negative
exponential moments of $X$, of orders $a$ and $-b$ respectively,
through the quantity
 $$ \left[\frac{{\cal
       L}(1+\delta)\beta)}{ 1-{\cal L}(-(1-\delta)\beta))
   }\right]^{\frac{1}{2\beta}}= \left[\frac{{\cal L}(a)}{ 1-{\cal
       L}(-b) }\right]^{\frac{1}{a+b}}.$$ 
 The random variable $T_R$ will intervene a random number of times,
 geometrically distributed with parameter approximately equal to
 $\eta^{R/L}$. The exponential rate in Theorem \ref{main} given by
$$\theta     \delta \beta =  \theta \frac{(a-b)}{2} $$ 
thus depends on the difference $(a-b)$, the Laplace transform
$\hat{\cal L}(u)=\left[\frac{{\cal L}(a)}{ 1-{\cal L}(-b)
  }\right]^{{\frac u 2} (1+\frac{c}{L})} \bar{{\cal L}}_{c+L}
\left(\frac{u (a+b)}{2} \right)$ of some r.v.  that accounts for the
time cost of failing a coupling attempt, and some small enough
$\theta\in(0,1)$ (so that $\hat{\cal L}(\theta) (1-\eta^{R/L})<1 $) that
controls the tradeoff between the previous ingredients.

We then deduce
 
\begin{coro}
  \label{corollary}
  For each $\gamma< \beta\theta\delta$ as in Theorem \ref{main}, there is an
  explicit constant $C$ depending on $\eta,L,c$ and $\gamma$ such that
  \[
  \| law(B_t^x) - law(B_t^0) \|_{TV} \leq \exp \left\{ \beta
    \theta x \right\} C \exp\left( -\gamma t
  \right).
  \]
  Moreover, if $U^x$ and $U^0$ respectively denote the renewal
  measures associated with the processes $(B_t^x)_{t\geq 0}$ and
  $(B_t^0)_{t\geq 0}$, then for all Borel sets $D$ in $(0,\infty)$ we
  have:
  \[
  \left| U^x(t+D)- U^0(t+D)\right|\leq 2 \exp \left\{ \beta \theta
    x \right\} C \exp\left( -\gamma t
  \right)(U^0((0,\sup D) )+1).
  \]
\end{coro}

\begin{rem}
  \label{refinedresults}
  By slight modifications of the proofs, it is also possible to
  replace the initial delay $0$ by a generic one $y\neq x$.  Moment
  conditions other than exponential can be treated with our techniques
  as well.
\end{rem}

In the next section, an outline of our approach and a plan of the
proofs are presented. A comparison to previous coupling arguments
together with a discussion of our results is given in {\bf Section 3}.

\section{Idea of the coupling and plan of the paper}

Our coupling construction and estimates will rely on the
discrete time structure of the renewal process. We start noting that,
under Assumption \ref{as:spread}, for any $s\in [0,L]$ and
$ t_2\geq t_1\geq s$ one has
\begin{equation*}
  \dP\big(X+s \in[t_1,t_2]\big)\geq \frac{\eta}{L}  
  \int_{t_1}^{t_2}\ind_{[c-L+s,c+L+s]}(u)\,du\geq  \frac{\eta}{L}  
  \int_{t_1}^{t_2}\ind_{[c,c+L]}(r)\,dr  .
\end{equation*}
The random variables $(X+s)_{s\in [0,L]}$ thus have a common uniform
component, of mass $\eta$, on the interval $[c,c+L]$. The following is
a straightforward and useful consequence:

\begin{lem}
  \label{comcomp}
  Under Assumption \ref{as:spread}, for each $s\in [0,L]$ one can
  define, on some probability space, a Bernoulli r.v.\ $\xi$ such
  that $\dP(\xi=1)=\eta=1- \dP(\xi=0)$, a uniform random variable $U$
  in $[c,c+L]$ independent of $\xi$, and random variables
  $X'\overset{d}{=}X$ and $X^{(s)}\overset{d}{=} X+s$, such that
  \begin{equation}
    \label{S'S''}
    \begin{split}
      \dP( X'=X^{(s)}  =U |\xi=1) & = 1 \, , \\
      \dP(X'\in dt | \xi=0) & = (1-\eta)^{-1} \left[\dP(X\in dt)
        - \eta \dP(U\in dt) \right]\, \forall \,   t\geq 0, \\
      \dP(X^{(s)} \in dt | \xi=0) & = (1-\eta)^{-1} \left[ \dP(X+s\in
        dt)- \eta\dP(U\in dt) \right] \, \forall \,  t\geq 0, \,
      \mbox{ with }\\
    \end{split}
  \end{equation}
  $X'$ and $ X^{(s)} $ independent conditionally on $\{\xi=0\}$. In
  particular, $(X',X'':=X^{(s)}-s)$ are two copies of the random
  variable $X$ for which $
  X''=X'-s $  holds  with probability $\eta.$
\end{lem}
Given a random variable $Z$, the same construction can be made
conditionally on $\{Z=s\}$, in which case Lemma \ref{comcomp} holds
true a.s.\ with respect to the law of $Z$, and $X'$ and $X''$ are
independent of $Z$ (though the pair $(X',X'')$ is not). Thus, starting
from a relative initial delay of $x>0$, by coupling
$k=\lceil x/L\rceil$ pairs of consecutive inter-arrivals of the two
processes, it is possible to produce simultaneous renewals with
probability at least $\eta^{\lceil x/L\rceil}>0$.  However, this
probability might be arbitrarily small if the initial relative delay
$x$ is not controlled, whereas, when such a ``coupling attempt''
fails, the resulting relative delay can in principle be arbitrarily
large.

Our coupling will therefore consist in a two-steps iterative scheme.
{\bf Step 1} roughly consists in running two independent copies until
renewals of both processes occur closer that some (large enough)
$R>0$.  The following bounds for the time this requires will be proved
in {\bf Section \ref{LyapRW}}:
  
\begin{prop}\label{lyap}
  Given $R>0$ and two independent copies ${(T_n)}_{n\geq 0}$ and
  ${(T'_n)}_{n\geq 0}$ of the renewal process such that $T_0=0$ and
  $T'_0=x>0$, let
  \begin{equation*}
    T_R=T_R(x):= \inf\{t\geq 0:  \exists n,n'\in \dN \mbox{ such that } 
    t=T_n \geq T'_{n'}-R \mbox{ or } t=T'_{n' } \geq T_n-R \}.
  \end{equation*} 
  Then, if
  $ R\geq \frac{1}{2\beta}\log \left[\frac{{\cal L}(\lambda+\beta)}{
      1-{\cal L}(-(\beta- \lambda)) }\right]$, for all
  $0<\lambda<\beta$ such that ${\cal L}(\lambda+ \beta)<\infty $, we
  have
  \begin{equation*}
    \dE_x(e^{\lambda T_R})\leq e^{\beta x \mathbf{1}_{x >R} } 
    \quad\mbox{ and }\quad  \forall \lambda ' \in [0,\lambda] , \, \,  
    \dE_x(e^{\lambda' T_R})\leq e^{\frac{\lambda' \beta}{\lambda} x
      \mathbf{1}_{x >R} } . 
  \end{equation*}
  If we moreover write $D_R(x)=|T_n-T'_{n'}|\leq R$, for $n,n'\in \dN$
  such that $T_R(x)=T_n\geq T'_{n'}-R$ or $T_R(x)=T'_{n'}\geq T_n-R,$
  the process $(T_R(x),D_R(x):x\geq 0)$ is measurable.
\end{prop}

Although the coupling in Step 1 is a classic one, the previous
exponential estimates are to our knowledge new; they rely on a
Lyapunov-type argument for some discrete-time random walk defined in
terms of the two copies' epochs.  Notice that Step 1 is not run (i.e.\
$T_R(x)=0$) if $x\leq R$.
 
As soon as the relative delay $z$ between the two copies is less than
$R$, {\bf Step 2} puts in place the coupling suggested after Lemma
\ref{comcomp}. More precisely, in Step 2 we will use the coupling of
two copies of the renewal process provided by the following result,
which is proved in {\bf Section \ref{lower bound}}:
  
\begin{lem}
  \label{momexpmax}
  For each $z>0$ one can define on some probability space two copies
  ${(T_n)}_{n\geq 0}$ and ${(T'_n)}_{n\geq 0}$ of the renewal process
  with $T_0=0$ and $T_0'=z$ and a random variable $I\in \dN$ a.s.
  bounded by $ \lceil z/L\rceil$, such that the event
  $ \Theta(z):=\left\{ T_{I}=T_{I}' \right\}$ and the random variable
  $M(z):=\max\{ T_{I},T_{I}' \} $ satisfy, for each $R>0$, the uniform
  bounds
  $$  \inf_{z\in [0,R]} \dP ( \Theta(z) ) \geq \eta^{ \lceil
    R/L\rceil}>0\, , $$     
  \begin{equation*}
    \sup_{z\in [0,R]}   \dE\left( e^{\gamma (M(z)- (z+ c \lfloor z/L\rfloor )  ) } \right)\leq   
    \bar{{\cal L}}_{c+L} (\gamma)
  \end{equation*}
  and
  \begin{equation*}
    \sup_{z\in [0,R]}     \dE\left( e^{\gamma (M(z)- (z+ c \lfloor z/L\rfloor )  ) }
      \ind_{\Theta(z)^c} \right)(1- \eta^{ \lceil z/L\rceil})^{-1}\leq    
    \bar{{\cal L}}_{c+L} (\gamma)
  \end{equation*}
  for all $\gamma \in \dR$.  Moreover, setting
  $m(z):=\min\{ T_{I},T_{I}' \} $, this construction can be done
  simultaneously for all $z>0$ in such a way that the process
  $(M(z),m(z),\ind_{\Theta(z)}:z\geq 0 )$ is measurable.
\end{lem}

The random variable $I$ in Lemma \ref{momexpmax} will correspond to
the minimal number of pairs of inter-arrivals, consecutively coupled
as in Lemma \ref{comcomp}, required to obtain simultaneous renewals
with positive probability, if $T_0=0$ and $T_0'=z$.  Thus, if Step 2
is run after Step 1, the event $\Theta(z)$ will occur with uniformly
lower bounded probability.  We say in that case that the coupling
succeeds; otherwise, one goes back to Step 1 and iterates.  
The upper bounds in Lemma \ref{momexpmax} moreover provide uniform
exponential estimates of the continuous time $M(z)$ spent during one
iteration of Step 2, whatever its outcome is, in terms of the relative
delay $z$ between the two copies at the beginning of it.  Thus, even
if their relative delay after Step 2 can be unbounded if the coupling
attempt fails, these bounds will provide some control of the initial
delay at the beginning of the next iteration of Step 1.

The proof of Theorem \ref{main}, given in {\bf Section
  \ref{proofmainresult}}, will consist in providing an exponential
control of the total continuous time required by the two copies,
constructed using this scheme, in order to have simultaneous renewals.
Hence, it will bring together Proposition \ref{lyap} and Lemma
\ref{momexpmax}, by means of an exponential estimate on
``sub-geometrical" sums of dependent positive random variables (Lemma
\ref{expdomgeom} in {\bf Appendix \ref{subgeomsum}}).  The first
statement of Corollary \ref{corollary} is straightforward from Theorem
\ref{main} and inequality \eqref{couplingineq}. The second one is more
subtle and is proved in {\bf Section \ref{bound renewal}}.
    
\section{Comparison to previous couplings and discussion of our
  results}

Proofs of renewal theorems given in \cite{LindvallRogers} or
\cite{alsmeyer}, among others, rely on the hitting times of intervals
$[-\varepsilon,\varepsilon]$ by a random walk $(T_n-T_n')_{n\geq 0}$
defined in terms of two renewal processes with coupled inter-arrivals
differing by less than $\varepsilon>0$. Those random walks being
symmetric, the expected number of steps in order that
$\varepsilon$-close renewals occur and the expected real time required
for that are infinite. Here, we will deal with a random walk which is
strongly biased towards $0$ and thus has some geometrically decaying
hitting times.  This walk is somehow reminiscent of a Markov chain
studied in \cite{Ney}, but our arguments are quite different and avoid
in particular the use of bounds based on the renewal equation.  Our
two-step coupling scheme is rather inspired by the celebrated
Meyn-Tweedie approach to long time convergence of Markov processes
(see \cite{MT}), but our discrete-time regenerative structure is
different. Our strategy also differs from the regenerative process
approach in continuous-time adopted in Ch.\ VII of \cite{asmussen03},
which at some point needs the use of asymptotic bounds on the renewal
measure and hence cannot yield tractable  estimates.

Another related reference is Chapter~6 in the book
\cite{kalashnikov94}. In this work, brought to the authors' attention
by an anonymous referee, Kalashnikov develops techniques close to the
ones of the present paper: he constructs a coupling under a condition
of contraction in total variation (condition (3) stated p.~167) and then
proves that this condition is satisfied when the inter-arrival times
have a distribution in some class (given in Definition 2, p.~185),
which is comparable to our Assumption~2. However, the author doesn't
follow all the constants and this turns out to be a difficult task (in
particular since Lemma~9 therein is given without proof). Thus, even if
the ideas in \cite{kalashnikov94} are close to ours, the present paper
provides a more direct approach (without any abstract contraction
condition), which makes it easier to exhibit bounds for the
convergence rate.

Unfortunately, the joint dependance of our bounds on the parameters
$( \beta, \theta, \delta)\in \dR_+\times (0,1)^2$ is not simple and,
in particular, the optimization problem one needs to solve in order
to maximize the convergence rate is not convex.  Although its
solution could be numerically approximated by some global optimization
routine, while optimizing the uniform component considered as well,  in general we
do not expect to get sharp bounds, since our arguments
rely on pessimistic (though careful) estimates. For instance, if $S$
has the folded standard Gaussian distribution, a numerical
optimization of our bounds yields the maximum rate
$\beta^* \theta^* \delta^*=0.001306$.  In turn, a nonlinear regression fit on
Monte-Carlo sample averages of quantities of the type $ U^x(t+D)$
suggest  in this case an exponential convergence rate about $3$ orders of
magnitude faster.

Nevertheless, the techniques here developed have the interest of
providing tractable bounds for the renewal theorem in a general
setting. In doing so they also give additional insight on the
properties of $S$ involved in the speed of convergence. Our arguments
could in principle be refined in order to take advantage of more
specific features of the inter-arrivals (such as the additional
integrability or increasing hazard rate of the above example). They
should allow for extensions to more general frameworks in renewal
theory as well.

%
%

\section{Proof of Theorem \ref{main}}\label{proofmainresult}

We start by estimating the total time spent during one iteration of
Step 2 followed by one of Step 1 (in that order), when at the
beginning of the former one of the two copies of the renewal process
is $0-$delayed and the other one has delay $z>0$.  In the notations of
Lemma \ref{momexpmax}, their relative delay at the end of Step 2 is
$\Delta(z):=M(z)-m(z)\geq 0, $ and one has $\Delta(z)=0$ if the
coupling is successful.  The total time spent in one iteration of
Steps 2 and then 1 has the same law as $T_R(\Delta(z))+ m(z),$ where
$x\mapsto T_R(x)$ is a copy of the (measurable) process of Proposition
\ref{lyap}, independent from the process
$z\mapsto (\Delta(z),m(z),\ind_{\Theta(z)})$.  For fixed
$0<\lambda'<\lambda<\beta$ and $R\geq 0$ as in Proposition \ref{lyap},
we get
\begin{equation*}
  \begin{split}
    \dE\left(\exp\{ \lambda'( T_R(\Delta(z))+ m(z))\}\right)=&\,
    \dE\left( \dE\left(\exp\{ \lambda'( T_R(x)\}
      \right)\bigg\vert_{x=\Delta(z)}
      \exp\{ \lambda' m(z)\} \right)\\
    \leq & \, \dE\left( \exp\left\{ \frac{\lambda' \beta}{\lambda}
        \Delta(z)+\lambda'  m(z) \right\} \right)\\
    \leq & \, \dE\left( \exp\left\{ \frac{\lambda' \beta}{\lambda}
        M(z)\right\} \right).\\
  \end{split}
\end{equation*}
We similarly obtain
\begin{equation*}
  \begin{split}
    \dE\left(\exp\{ \lambda'( T_R(\Delta(z))+
      m(z))\}\ind_{\Theta(z)^c}\right)=&\, \dE\left( \dE\left(\exp\{
        \lambda'( T_R(x)\} \right)\bigg\vert_{x=\Delta(z)}
      \exp\{ \lambda' m(z)\}\ind_{\Theta(z)^c} \right)\\
    \leq & \, \dE\left( \exp\left\{ \frac{\lambda' \beta}{\lambda}
        M(z)\right\}\ind_{\Theta(z)^c} \right).\\
  \end{split}
\end{equation*}
Now,  by Lemma \ref{momexpmax},
$\inf_{z\in [0,R]} \dP ( \Theta(z) ) \geq \eta^{ \lceil R/L\rceil} $
 for each $R>0$. Taking
$\gamma = \frac{\lambda' \beta}{\lambda}$ therein we get
\begin{equation}
  \label{estimexpoM}
  \sup_{z\in [0,R]} \dE\left( \exp\{ \lambda'( T_R(\Delta(z))+ m(z))\}\right)
  \leq   e^{\frac{\lambda' \beta}{\lambda} (R +\lfloor R/L\rfloor c)}   
  \bar{{\cal L}}_{c+L} \left(\frac{\lambda' \beta}{\lambda}\right)
\end{equation}
and
\begin{equation}\label{estimexpoMfail}
  \sup_{z\in [0,R]}  \dE\left( \exp\{ \lambda'( T_R(\Delta(z))+ m(z))\}
    \ind_{\Theta(z)^c} \right)\leq  e^{\frac{\lambda' \beta}{\lambda} (R +\lfloor R/L\rfloor c)}   
  \bar{{\cal L}}_{c+L} \left(\frac{\lambda' \beta}{\lambda}\right)
  (1-\eta^{ \lceil R/L\rceil}).
\end{equation}

Let us now derive an exponential estimate for the global time required
for the two copies in our coupling scheme to have simultaneous
renewals. Thanks to the independence of the inter-arrivals of the
renewal process and the measurability properties stated in Proposition
\ref{lyap} and Lemma \ref{momexpmax}, the relevant time-lengths in our
coupling scheme can be constructed using independent sequences
$(( T_R^j(y),D_R^j(y)): y\geq 0)_{j\in \dN}$ and
$( (\Delta_j(z),m_j(z),\ind_{\Theta_j(z)}): z\geq 0)_{ j\in \dN
  \backslash \{0\}}$
of independent copies of the processes $(T_R,D_R)$ and
$(\Delta,m,\ind_{\Theta})$. More precisely, recursively defining
\[
Y_0=x, \quad Z_{j+1}=D^j_R(Y_j) \mbox{ and } Y_{j+1}=\Delta_{j+1}
(Z_{j+1}), \quad j\geq 0,
\]
the sequence $(Y_j,Z_{j+1})_{ j\in \dN}$ has the same law as the
sequence of relative delays of the two copies, after the $j-$th
iteration of Step 1 and the consecutive one of Step 2, respectively. A
stochastic upper bound for the coalescing time of the two copies is
then given by
\[
\bar{T}^*(x):=T_R^0(x)+\sum_{j=1}^{\sigma}
\left(T_R^j(Y_j)+m_j(Z_j)\right) \, ,
\]
where $\sigma=\inf\{j>0 \, : Y_j=0 \}$.  Applying Lemma
\ref{expdomgeom} in Appendix \ref{subgeomsum} to the filtration
$({\cal G}_n)_{n\in \dN}$, with 
\[ {\cal G}_n:=\sigma\left(\left\{ T^j_R, D_R^j,
    \Delta_k,m_k,\ind_{\Theta_k} : j =0,\dots n \, , k=1,\dots, n
  \right\}\right)
\] 
and the random variables and events $W_n= T_R^n(Y_n)+m_n(Z_n)$ and
$A_n=\{ \ind_{\Theta_n(Z_n)}=1 \}$, we deduce, thanks to independence
of the processes generating ${\cal G}_n $ and the bounds
\eqref{estimexpoM} and \eqref{estimexpoMfail}, that
\begin{equation*}
  \begin{split}
    \dE\left( e^{\lambda' \bar{T}^*(x)}\right)\leq & \, \dE\left(
      e^{\lambda' T_R^0(x)}\right) \dE\left(\left(
        e^{\left\{\frac{\lambda' \beta}{\lambda} (R +\lfloor
            R/L\rfloor c)\right\}} \bar{{\cal L}}_{c+L}
        \left(\frac{\lambda' \beta}{\lambda}
        \right)  \right)^G \right) \\
    \leq & \, e^{\frac{\lambda' \beta}{\lambda} x \mathbf{1}_{x >R} }
    \dE\left( \left( e^{\left\{\frac{\lambda' \beta}{\lambda} (R
            +\lfloor R/L\rfloor c)\right\}} \bar{{\cal L}}_{c+L}
        \left(\frac{\lambda' \beta}{\lambda}\right)
      \right)^G \right)  \\
  \end{split}
\end{equation*}
where $G$ is a geometric r.v.\ of parameter
$\eta^{ \lceil R/L\rceil}\in (0,1)$.  Given parameters $\beta,\theta$
and $\delta$ as in Theorem \ref{main}, its proof  is
then achieved by taking above $\lambda=\delta \beta$,
$\lambda'=\theta \delta \beta$ and $ R=R(\delta,\beta)$.

\section{Step 1: a positive recurrent random walk associated with
  independent renewal processes}\label{LyapRW}

We next prove Proposition \ref{lyap}.  We introduce to that end a
biased random walk $(Y_n)_{n\in \dN}$ in $\dR$ defined from a single
sequence of i.i.d.\ inter-arrivals $(X_n)_{n\in \dN}$. The process
$(Y_n)_{n\in \dN}$ will account for the relative signed (positive or
negative) delay of one fixed copy of the renewal process with respect
to a second copy, after a total number $n$ of inter-arrivals has
occurred. More precisely, given an initial relative delay $x\in \dR$,
we set $Y_0=x$.  By convention, $x\geq 0$ means that one copy,
henceforth fixed and called ``the first copy'', has a delay $x$,
whereas the other copy, called ``the second copy'', is $0$ delayed.
Conversely, a relative initial delay $x<0$ means that the first copy
is $0$ delayed and the second one has a delay $|x|>0$.  To construct
the walk we proceed as follows: if for given $n$ we have
$Y_{n}\geq 0$, meaning that the first copy's last defined epoch
occurred at distance $Y_{n}$ to the right of the second copy's one, we
add the next inter-arrival $X_{n+1}$ to the last defined epoch of the
second copy. If, on the contrary, we had $Y_{n}< 0$, this means that
the first copy's last defined epoch occurred at distance $|Y_{n}|$
left from the second copy's one, and the random variable $X_{n+1}$ is
then added to the last defined epoch of the first copy.  We thus have
\[
Y_{n+1}=
\begin{cases}
  Y_n-X_{n+1} & \text{if $Y_n\geq  0$,}\\
  Y_n+X_{n+1} & \text{if $Y_n<0.$}\\
\end{cases}
\]
Notice that the ``leftmost copy'' by the end of step $n$ either
catches up in step $n+1$ part of its delay with respect to the other
copy or otherwise overshoots the lastly defined epoch of the latter,
in which case the roles are then interchanged.  Setting
$N_n^+:=\inf\{ m\in \dN: \sum_{i=0}^m \ind_{Y_{i-1}\geq 0} \geq n\} $
and
$N_n^-: =\inf\{ m\in \dN: \sum_{i=0}^m \ind_{Y_{i-1}< 0}\geq n \} $,
it easily follows from the independence of the $(X_n)$ that $N_n^+$
and $N_n^-$ go to $\infty$ with $n$. Moreover, the inter-arrivals
assigned to the first and second copies are respectively given by the
sequences $ (X_{N^+_n})_{n\geq 1} \mbox{ and }(X_{N^-_n})_{n\geq 1}, $
and the strong Markov property of the random walk $(Y_n)_{n\in \dN}$
shows that these are independent i.i.d.\ sequences; they thus define
independent copies of the renewal process.

The minimal total number of inter-arrivals required for epochs of
these two processes to take place not farther that $R$ from each
other is
$$\tau_R:=\inf\{ n\geq 0: |Y_n|\leq R\}.$$
Hence, in the notation of Proposition \ref{lyap}, we have
\begin{equation*}
  T_R(x)= \mathbf{1}_{|x| >R}\left(  \sum_{i=1}^{\tau_R} X_i\ind_{Y_{i-1}\geq 0} 
    +  Y_{\tau_R}\ind_{Y_{\tau_R}< 0}   \right) \, \leq  \bar{T}_R(x):=
  \mathbf{1}_{|x| >R} \sum_{i=1}^{\tau_R} X_i\ind_{Y_{i-1}\geq0},
\end{equation*}
and $D_R(x)=|Y_{\tau_R}|$. We will estimate exponential moments of
$\bar{T}_R(x)$ by a Lyapunov-type argument.  Let $0\leq \lambda<\beta$
be such that ${\cal L}(\lambda+ \beta)<\infty$ and set
$V(x):=e^{\beta |x|}$. Then,
\begin{equation*}
  \begin{split}
    \dE_x\Big( V(Y_1)& e^{ \lambda X_1\ind_{x\geq 0}}\Big)\\
    \leq & e^{-\beta |x|}\left[\ind_{x\geq 0} \dE(e^{(\lambda
        +\beta)X}) + \ind_{x< 0} \dE(e^{\beta X}) \right]+ e^{\beta
      |x|}\left[\ind_{x< 0} \dE(e^{-\beta X})
      + \ind_{x\geq  0} \dE(e^{(\lambda -\beta) X}) \right]\\
    \leq & e^{-\beta |x|} \dE(e^{(\lambda +\beta)X}) +
    e^{\beta |x|}  \dE(e^{(\lambda -\beta) X})\\
    \leq &V(x)\Big({\cal L}(-(\beta- \lambda) ) +e^{-2\beta R}{\cal
      L}(\lambda+\beta)\Big)+
    {\cal L}(\lambda+\beta)\ind_{[0,R]}(|x|)\\
  \end{split}
\end{equation*}
where the first inequality is obtained after partitioning the
expectation according to the signs of $x-X_1$ and of $x$.  By standard
arguments, the above bound entails that the discrete time process
\[
V( Y_{\tau_R\wedge n } ) e^{ \lambda \sum_{i=1}^{\tau_R\wedge n }
  X_i\ind_{\{Y_{i-1}\geq 0\}}}\rho^{- \tau_R\wedge n
}_{\beta,\lambda,R} \, , \quad n \in \dN,
\] 
with
\[
\rho_{\beta,\lambda,R}:= {\cal L}(-(\beta- \lambda) )+e^{-2\beta
  R}{\cal L}(\lambda+\beta),
\] 
is a supermartingale in the discrete filtration generated by the
sequence $(X_i)_{i\geq 1}$.  In particular,
\[
\dE_x\left( e^{ \lambda \sum_{i=1}^{\tau_R\wedge n }
    X_i\ind_{Y_{i-1}\geq 0}} \rho^{- \tau_R\wedge n
  }_{\beta,\lambda,R}\right)\leq \dE_x\left(V( Y_{\tau_R\wedge n } )
  e^{ \lambda \sum_{i=1}^{\tau_R\wedge n } X_i\ind_{Y_{i-1}\geq 0}}
  \rho^{- \tau_R\wedge n }_{\beta,\lambda,R}\right)\leq V(x).
\]
By letting $n\to \infty$ in the first expectation above, we deduce
that $\tau_R<\infty$ a.s.\ if $\rho^{- 1}_{\beta,\lambda,R}>1$ or if
$\rho^{- 1 }_{\beta,\lambda,R}=1$ and $0<\lambda<\beta$ (in the second
case we use the fact that
$\sum_{i=1}^{n} X_i \ind_{Y_{i-1}\geq 0} =\sum_{i=1}^{n^+ } X_{N_i^+}$
is a sum of i.i.d.\ random variables). This yields
$$\dE_x\left( e^{ \lambda \bar{T}_R(x) }\right)\leq e^{\beta x
  \mathbf{1}_{x >R} }$$ whenever $\rho_{\beta,\lambda,R} \leq 1$ or,
equivalently, when
$ R\geq\frac{1}{2\beta}\log \left[\frac{{\cal L}(\lambda+\beta)}{
    1-{\cal L}(-(\beta- \lambda) )}\right]=R(\lambda/\beta,\beta).  $
The first assertion of Proposition \ref{lyap} follows. The second one
is easily obtained with Holder's inequality.  The last assertion of
Proposition \ref{lyap} is straightforward from the previous
construction.


\section{Step 2: attempting an exact coupling  }\label{lower bound} 

We first construct the coupling of Lemma \ref{momexpmax}, in such a way
that the measurability condition in its last assertion is granted from
the beginning; we then establish the claimed exponential estimates.

Consider the four independent i.i.d.\ sequences:
$(\xi_i)_{i=1}^{\infty}$ of Bernoulli r.v.\  with
$\dP(\xi_i=1)=\eta = 1-\dP(\xi_i=0)$, $(U_i)_{i=1}^{\infty}$ of
uniform r.v.\ in $[c,c+L]$, and $(W'_i)_{i=1}^{\infty}$ and
$(\hat{W}_i)_{i=1}^{\infty}$ of r.v.\ such that
\begin{equation*}
  \begin{split}
    \dP(W'_i\in dt ) & = (1-\eta)^{-1} \left[ \dP(X\in dt)
      - \eta \dP(U\in dt)\right] ,   \,  t\geq 0   \, \mbox{ and } \\
    \dP(\hat{W}_i \in dt ) & = (1-\eta)^{-1} \left[
      \dP(X+L\in dt)-\eta \dP(U\in dt)\right]    , \,  t\geq 0, \\
  \end{split}
\end{equation*}
with $U$ uniformly distributed in $[c,c+L]$. Consider also $\vartheta$
a uniform random variable in $[0,1]$ independent of all the previous
ones. By Lemma 3.22 in \cite{Kall} there exists a measurable function
$\Phi: \dR_+\times [0,1]\to \dR_+$ such that, for each $z\in\dR_+$,
the random variable $\Phi(z,\vartheta)$ satisfies
\begin{equation*}
  \dP(\Phi(z,\vartheta)  \in  dt  )= (1-\eta)^{-1}  \left[  
    \dP(X+(z- L \lfloor z/L\rfloor) \in dt)-\eta \dP(U\in dt)\right]  
  \,  ,  t\geq 0 .
\end{equation*}
Set  now $k= k(z):=\lceil z/L \rceil$ and for 
$i\geq 1$ define:
\begin{equation*}
  \begin{split}
    X'_i := & \, \ind_{\xi_i=1} U_i +  \ind_{\xi_i=0} W'_i \, ,  \\
    \hat{X}_i := & \, \ind_{\xi_i=1} U_i + \ind_{\xi_i=0} \left(
      \ind_{i< k(z)}
      \hat{ W}_i + \ind_{i= k(z)}  \Phi(z,\vartheta) \right)   \mbox{   and}\\ 
X''_i:=  & \hat{X}_i - \left( \ind_{i< k(z)} L + \ind_{i= k(z)} (z- L
  \lfloor z/L\rfloor) \right) .
\end{split}
\end{equation*}
We remark for later use that the r.v.\ $W_i''$ defined as
$W''_i:=\hat{W}_i-L$, for $1\leq i<k$, and as
$W''_k:=\Phi(z,\vartheta)-(z-L \lfloor z/L\rfloor)$, all have the same
law as the r.v.\ $W'_i$.

The sequences $( X_j')_{j=1}^k$ and $( X_j'')_{j=1}^k$ are both
i.i.d.\ with the same law as $X$, they are measurable functions
jointly in $z$ and randomness and, on 
$F_k:=\{(\xi_1,\dots,\xi_k)=(1,\dots,1)\}$,  it a.s.\ holds that
\begin{equation*}\label{S'S''Z}
  X'_1+ \dots+X'_{k}-( X''_1+\dots+X''_k) = z.
\end{equation*}
In particular, the probability of having such an equality is bounded
from below by $\eta^{ \lceil z/L\rceil}>0.$ The coupling will then
consist in sampling the random variables $\xi_i, X'_i$ and $X''_i$
up to the random index
$$I:=\inf\{ j\geq 0 : \xi_j=0\}\wedge k\leq \lceil z/L \rceil.$$ 
If the latter set is empty, then the event $F_k$ occurs, simultaneous
renewals take place at time $X'_1+ \dots+X'_{k}$ and the coupling
attempt is successful; otherwise, we say that it fails. Notice that if
$1\leq I<k$, the coupling attempt is said to fail, even if the $I$-th
renewals of the two copies take place simultaneously (which can for
instance happen if $X$ has atoms). Notice also that when $I=k$, the
coupling might succeed or fail. In all cases, we have
$$ M(z)=\max\left\{\sum_{j=1}^{I} X'_j, z+\sum_{j=1}^{I} X''_j\right\}.  $$
The random variable $m(z)$ in the statement corresponds to
$m(z)=\min\left\{\sum_{j=1}^{I} X'_j, z+\sum_{j=1}^{I} X''_j
\right\}$, and the event $ \Theta(z)$, which corresponds to
$\left \{ M(z)=m(z)\right\} $, occurs if $F_k$ does.  The indicator
function in the second estimate in Lemma \ref{momexpmax} can thus be
replaced by $\ind_{F_k^c}$.  Since $X''_j\leq c$ for $j\leq I$ and
$I\leq \lfloor z/L\rfloor$, we always have (with $\sum_{\emptyset}=0$)
\[
\sum_{j=1}^{I-1} X'_j\leq z+\sum_{j=1}^{I-1} X''_j \leq z+ c
\lfloor z/L\rfloor.
\]
 
Moreover, on $F_k$ we have $X'_{I}, X''_{I} \leq (c+L)$. It then
follows on one hand that, for all $\gamma\in \dR,$
\begin{equation*}
  \begin{split}
    \dE\left( e^{\gamma M(z)} \right)\leq & e^{\gamma (z+ c \lfloor
      z/L\rfloor ) } \left[e^{\gamma (c+L)} \dP(F_k) +\dE\left(
        e^{\gamma\max\{ X'_{I}, X''_{I} \} }
        \vert F_k^c\right)\dP(F_k^c)  \right]\\
    \leq & e^{\gamma (z+ c \lfloor z/L\rfloor ) } \max\left\{
      e^{\gamma (c+L)} ,
      \dE\left( e^{\gamma\max\{ X'_{I}, X''_{I} \} } \vert F_k^c\right)\right\} \\
    \leq & e^{\gamma(z+ c \lfloor z/L\rfloor ) }
    \dE\left( e^{\gamma\max\{  (c+L), X'_{I}, X''_{I} \} } \vert F_k^c\right)\leq \infty.  \\
  \end{split}
\end{equation*}
On the other hand, we obtain
\begin{equation*}
  \begin{split}
    \dE\left( e^{\gamma M(z)} \ind_{F_k^c}\right)\leq & e^{\gamma(z+ c
      \lfloor z/L\rfloor ) }
    \dE\left( e^{\gamma\max\{ X'_{I}, X''_{I} \} } \vert F_k^c\right)\dP(F_k^c)  \\
    \leq & e^{\gamma(z+ c \lfloor z/L\rfloor ) } \dE\left(
      e^{\gamma\max\{ (c+L), X'_{I}, X''_{I} \} } \vert
      F_k^c\right)\dP(F_k^c)
    \leq \infty.  \\
  \end{split}
\end{equation*}
The two required estimates will then be proved by showing that
\begin{equation}\label{SSWW}
  \dE\left( e^{\gamma\max\{  (c+L), X'_{I}, X''_{I} \} } \vert F_k^c\right)\leq  
  \dE\left(  e^{\gamma\max\{   W', W'' \}} \right), 
\end{equation}
for independent r.v.  $( W', W'')$ of law
$ \left(\frac{\dP(X\leq c+L)-\eta}{1-\eta}\right)\delta_{ c+L}(ds)+
\frac{ \ind_{ s>c+L}}{1-\eta}\dP_X(ds) $ with $\dP_X$ the law of $X$.
Indeed, since $ \dP(X>c+L)\leq 1-\eta$, one gets
$\dP(W'>s)\leq \dP(X>s\vert X>c+L)$ for all $s\geq 0$, that is, $W'$
is stochastically smaller than a r.v.  $\bar{X}'$ equal in law to $X$
conditioned on being not smaller than $c+L$.  It then follows that
$\max\{ W', W'' \}$ is stochastically smaller than
$ \max\{ \bar{X}', \bar{X}'' \} $ for an i.i.d.  pair
$(\bar{X}', \bar{X}'')$, from where we conclude.

Let us thus check inequality \eqref{SSWW}.  Since
$ \ind_{F_k^c}=\sum_{l=1}^{k-1} \ind_{I=l} + \ind_{I=k,\xi_k=0} $,
we have
\begin{equation*}
  \begin{split}
    \dE\left( e^{\gamma\max\{ (c+L), X'_{I}, X''_{I} \} }
      \ind_{F_k^c}\right)= & \sum_{l=1}^{k-1} \dP(I=l) \dE\left(
      e^{\gamma\max\{ (c+L), W'_{l}, W''_{l} \} }
    \right) \\
    & + \dP(I=k,\xi_k=0) \dE\left( e^{\gamma\max\{ (c+L), W'_k,
        W''_k \}}
    \right)   
  \end{split}
\end{equation*}
so it suffices to show that, for $l=1,\dots, k$,
$\dE\left( e^{\gamma\max\{ (c+L), W'_l, W''_l \}} \right)$ is bounded
by $ \dE\left( e^{\gamma\max\{ W', W'' \}} \right) $.  This follows
from
\begin{equation*}
  \begin{split}
    \dE\left( e^{\gamma\max\{ (c+L), W'_l, W''_l \}} \right) = &
    e^{\gamma (c+L)} \dE\left( \ind_{ W'_l, W''_l\leq c+L}
    \right)+ \dE\left( e^{\gamma W'_l}
      \ind_{ W'_l>c+L\geq  W''_l }  \right) \\
    &+ \dE\left( e^{\gamma W''_l} \ind_{ W''_l>c+L\geq W'_l }
    \right) + \dE\left( e^{\gamma\max\{ W'_l, W''_l \}}
      \ind_{ \min\{W'_l, W''_l \}>c+L}   \right).  \\
  \end{split}
\end{equation*}
and the fact that, for each $l=1,\dots k-1$, $(W'_l, W''_l )$ are
independent,
$\dP(W'_l\leq c+L)=\dP( W''_l\leq c+L)=\frac{\dP(X\leq
  c+L)-\eta}{1-\eta}$
and
$\dE( f(W'_l) \ind_{ W'_l>c+L})=\dE( f( W''_l) \ind_{
  W''_l>c+L})=\frac{ \dE( f(X) \ind_{ X>c+L})}{1-\eta}$
for all nonnegative measurable function $f$.

\section{Bounds for the renewal measure}\label{bound renewal}

Thanks to inequality \eqref{couplingrenew} and the fact that
$t+D\subseteq (t,t+h] $ for $h=\sup D$, to prove the second statement
of Corollary \ref{corollary} it is enough to show that, for any $h>0$,
\begin{equation}
  \label{boundsumrenew} 
  \dE\left(\ind_{T^*(x)>t}  \left(\sum_{j=0}^{\infty}\ind_{(t,t+h]}(T''_j)\right)\right) 
  \leq \dP (T^*(x)>t)  (U^0((0,h])+1) 
\end{equation}
for $(T''_n)=(T_n)$ and $(T''_n)=(T'_n)$ the epochs of the two
copies. To that end we describe the discrete time structure used in
constructing our coupling in a slightly different way from before.
Consider the following independent i.i.d.\ sequences:
\begin{itemize}
\item[$\bullet$] $(\tilde{X}_k)_{k=1}^{\infty}$ with law equal to that
  of $X$,
\item[$\bullet$] $(U_k)_{k=1}^{\infty}$ uniform in $[c,c+L]$,
\item[$\bullet$] $(\xi_k)_{k=1}^{\infty}$ Bernoulli of parameter
  $\eta$,
\item[$\bullet$] $(W'_k)_{k=1}^{\infty}$ and
  $(\hat{W}_k)_{k=1}^{\infty}$ with the laws described in Section
  \ref{lower bound} and
\item[$\bullet$] $(\vartheta_k)_{k=1}^{\infty}$ uniform in $[0,1]$.
\end{itemize}

We can then construct our coupling using the i.i.d.\ random vectors
$(\tilde{X}_k, U_k, \xi_k, W'_k,\hat{W}_k, \vartheta_k)_{k\in \dN}$,
as follows. We run Step 1 using the random variables $(\tilde{X}_k)$
to construct the random walk of Section \ref{LyapRW}, until the
conditions required to start Step 2 (i.e. a relative delay not larger
that $R$) are met. This first happens at some discrete random time,
which is a stopping time with respect to the filtration
$({\cal F}_m)_{m\geq 1}$ defined by
$${\cal F}_m:=\sigma \left(  \tilde{X}_k, U_k, \xi_k,  W'_k,\hat{W}_k, 
  \vartheta_k: k=1,\dots, m \right).$$
Notice that, until then, the remaining coordinates
$( U_k, \xi_k, W'_k,\hat{W}_k, \vartheta_k)$ of the vector are not
used. Moreover, one and only one copy of the renewal process has had a
renewal at each time step $k$.  Right after that stopping time, we
start Step 2 using at each time step $k$ some random variables among
$U_k, \xi_k, W'_k,\hat{W}_k $ and $\vartheta_k$ (as needed in the
scheme described in Section \ref{lower bound}).  This is done until
some second stopping time (with respect to $({\cal F}_m)_{m\geq 1}$)
at which the coupling attempt succeeds or fails. In the latter case
one restarts Step 1. Notice that during Step 2, both copies have one
renewal at each time step $k$.

We denote by $\tau(1)<\tau(2)<\tau(3)<\cdots$ (resp.\
$\tau'(1)<\tau'(2)<\tau'(3)<\cdots$) the discrete times $k$ at which
the number of arrivals of the first (resp.\ second) copy of the
renewal process is increased by one additional unit.  Notice they are
also stopping times with respect to $({\cal F}_m)_{m\geq 1}$.

We then denote by $X_n$ (resp.\ $X_n'$) the increment of the first
(resp.\ second) copy at time $k=\tau(n)$ (resp.\ $k=\tau'(n)$).  It is
then not hard to see that the sequence $(X_n, X_n')_{n\geq 1} $ has the same
law as the sequence of pairs  of  inter-arrivals resulting from  our coupling
construction.  Moreover, $(X_n)_{n\geq 1} $ and $(X_n')_{n\geq 1} $
are respectively adapted to the filtrations
$({\cal F}_{\tau(n)})_{n\geq 1}$ and
$({\cal F}_{\tau'(n)})_{n\geq 1}$. 

Denote by $(T_n)$ and $N_t$ the epochs and counting processes
corresponding to this sequence $(X_n)$ and observe that
$\{N_t=n\} =\{T_n\leq t \} \cap \{T_{n+1}>t \}\in {\cal
  F}_{\tau(n+1)}$. Moreover,
$\{T^*(x)>t, N_t=n\}\in {\cal F}_{\tau(n+1)}$ since this event can be
written in terms of $\{N_t=n\} $ and the family of random variables
$\left(\ind_{\{k\leq \tau(n+1)\}} ( \tilde{X}_k, U_k, \xi_k,
  W'_k,\hat{W}_k, \vartheta_k)\right)_{k\geq 1}$

Defining now a function $F$ on $[0,\infty)^{\dN\backslash \{0\}}$ by
$F(x_1,x_2, \dots )= \sum_{j=1}^{\infty}\ind_{(0,h]}(\sum_{k=1}^j
(x_k))$, the expectation in the left hand side of
\eqref{boundsumrenew} is seen to be equal to
\begin{equation*}
  \begin{split}
    \sum_{n\in \dN} \dE\Bigg(\ind_{\{T^*(x)>t, N_t=n\}}&
    \Big(\sum_{j=n+1}^{\infty}\ind_{(t,t+h]}(T_j)\Big) \Bigg)\\
    = & \sum_{n\in \dN} \dE\left(\ind_{\{T^*(x)>t, N_t=n\}}
      F(T_{n+1}-t, X_{n+2},X_{n+3},\dots  ) \right) \\
    \leq & \sum_{n\in \dN} \dE\left(\ind_{\{T^*(x)>t, N_t=n\}}
      (1+F( X_{n+2},X_{n+3},\dots  )) \right). \\
\end{split}
\end{equation*}

To conclude \eqref{boundsumrenew} for $(T''_n)=(T_n)$ it suffices to check that
$ \dE\left(\ind_{\{T^*(x)>t, N_t=n\}}F( X_{n+2},X_{n+3},\dots )
\right)= \dP \{T^*(x)>t, N_t=n\} \dE\left( F(X_1,X_2, ...) \right)$
for all $n\in \dN$. This property is a consequence of the strong
Markov property of the (i.i.d.) process
$\left( \tilde{X}_k, U_k, \xi_k, W'_k,\hat{W}_k,
  \vartheta_k\right)_{k\geq 1}$ since, for each $n\in \dN$,
$\{T^*(x)>t, N_t=n\}\in {\cal F}_{\tau(n+1)}$ and the r.v.
$X_{n+2},X_{n+3},\dots$ can be constructed using the random vectors
$\left(\ind_{\{k\geq \tau(n+2)\}} ( \tilde{X}_k, U_k, \xi_k,
  W'_k,\hat{W}_k, \vartheta_k)\right)_{k\geq 1}$. The proof for
$(T''_n)=(T'_n)$ is similar.

  \appendix

\section{Appendix}
       
\subsection{ Laplace bounds for sub-geometric sums of dependent random
  variables }
\label{subgeomsum}
      
\begin{lem}
  \label{expdomgeom} 
  Let ${(\cal G}_n)_{n\geq 0}$ be a filtration, and $(A_n)_{n\geq 1}$
  and $(W_n)_{n\geq 1}$ sequences of respectively adapted events and
  adapted nonnegative random variables. Let
  $\sigma:=\inf\{ n\geq 1: \, \mathbf{1}_{A_n}=1\}$ and assume there
  exist $p\in (0,1)$ and a function $\psi:I\to \dR$ defined in some
  real interval $I$, such that for all $\lambda \in I$,
\begin{itemize}
\item[i)]
  $ \dE\left( e^{\lambda W_n} \vert { \cal G}_{n-1}\right) \leq
  e^{\psi(\lambda)} $
  a.s. on $A_1^c\cap\cdots \cap A_{n-1}^c$ if $n>1$ and a.s.\ if $n=1$
  and
\item[ii)]
  $ \dE\left( e^{\lambda W_n} \mathbf{1}_{A_n^c} \vert {\cal
      G}_{n-1}\right) \leq (1-p) e^{\psi(\lambda)} $
  a.s.\ on $A_1^c\cap\cdots \cap A_{n-1}^c$ if $n>1$.
\end{itemize}
Then,
$ \dE\left( e^{\lambda \sum_{j=1}^\sigma W_j} \right) \leq \dE\left(
  \left( e^{\psi(\lambda)} \right)^G \right)$
for all $\lambda \in I$, where $G$ is geometric of parameter $p$.  In
particular, if $0\in I$ and $\psi$ is increasing and goes to $0$ at
$0$, we have
$ \dE\left( e^{\lambda \sum_{j=1}^\sigma W_j} \right)\leq \frac{p
  e^{\psi(\lambda)}}{1- e^{\psi(\lambda)}(1-p)}$
for all $\lambda$ such that $\psi(\lambda)<- \log (1-p)$.
\end{lem}

If $\psi(0)=0$, condition ii) classically yields that $\sigma$ is
stochastically smaller than $G$ (see e.g.\ Lemma A.6 in
\cite{asmussen03}). The bound in Lemma \ref{expdomgeom} is sharp given
the assumptions (it is attained for $(W_n)$ i.i.d.\ of exponential law
with Laplace transform $e^{\psi(\lambda)}$ and $\sigma=G$
independent).

\dem We may assume that $\psi(\lambda)$ is
in the domain of the Laplace transform of $G$. Moreover, replacing
$(A_n)$ by $(\tilde{A}_n)$ defined as $\tilde{A}_n=A_n$ if $n\leq N$
and $\tilde{A}_n=\Omega$ if $n\geq N+1$, for some fixed integer $N$,
we may assume that $\sigma =\sigma\wedge N \leq N$ and then pass to
the general case using monotone convergence.  The fact that $\sigma$
is bounded justifies the interchange of sums with differences needed
to get (with the convention $\sum_{j=1}^0 = 0$) :
\begin{align*}
  \dE\left( e^{\lambda \sum_{j=1}^\sigma W_j} - 1 \right)  = & \sum_{n=1}^{\infty} 
  \dE\left( \left[ e^{\lambda \sum_{j=1}^n  W_j} -e^{\lambda \sum_{j=1}^{n-1}  W_j}\right]   
    \mathbf{1}_{\sigma \geq n} \right) \\ = & \sum_{n=1}^{\infty} 
  \dE\left( e^{\lambda \sum_{j=1}^{n-1}  W_j} \dE\left(e^{\lambda W_n}-1\vert {\cal G}_{n-1}\right) 
    \mathbf{1}_{\sigma \geq n} \right).
\end{align*}
We deduce, using i) to get the first inequality and ii) to get the
second one, that
\begin{equation*}
  \begin{split}
    \dE\Big( & e^{\lambda \sum_{j=1}^\sigma W_j} - 1 \Big)\\ & \leq
    \sum_{n=1}^{\infty} \left[e^{\psi(\lambda)} -1\right] \dE\left(
      \mathbf{1}_{n=1} + e^{\lambda \sum_{j=1}^{n-1} W_j}
      \mathbf{1}_{A_{n-1}^c} \cdots   \mathbf{1}_{A_{1}^c} \mathbf{1}_{n\geq 2} \right) \\
    & = \sum_{n=1}^{\infty} \left[e^{\psi(\lambda)} -1\right]
    \dE\left( \mathbf{1}_{n=1} + \dE\left( e^{\lambda W_{n-1}}
        \mathbf{1}_{A_{n-1}^c} \vert {\cal G}_{n-2}\right)
      \left(\mathbf{1}_{n=2} + e^{\lambda \sum_{j=1}^{n-2} W_j}
        \mathbf{1}_{A_{n-2}^c} \cdots
        \mathbf{1}_{A_{1}^c}\mathbf{1}_{n\geq 3}    \right)  \right) \\
    & \leq \sum_{n=1}^{\infty} \left[e^{\psi(\lambda)} -1\right]
    \mathbf{1}_{n=1} + \left[e^{2 \psi(\lambda)} -e^{\psi(\lambda)}
    \right] \left( 1-p \right) \dE\left( \mathbf{1}_{n=2} + e^{\lambda
        \sum_{j=1}^{n-2} W_j}
      \mathbf{1}_{A_{n-2}^c} \cdots   \mathbf{1}_{A_{1}^c}\mathbf{1}_{n\geq 3}   \right).
  \end{split}
\end{equation*}
Conditioning on $ {\cal G}_{n-3}$ in the last expectation and
iterating the argument yields the upper bound
\begin{equation*}
  \sum_{n=1}^{\infty} 
  \left[ e^{ n   \psi(\lambda) } -e^{ (n-1)   \psi(\lambda)  }\right] (1-p)^{n-1}   =   \sum_{n=1}^{\infty} 
  \left[ e^{ n   \psi(\lambda) } -e^{ (n-1)   \psi(\lambda)  }\right] \dP(G\geq n)  =  \dE\left(  \left( e^{\psi(\lambda)}  \right)^G- 1 \right).
\end{equation*}

{\bf Acknowledgements:} We thank anonymous referees for valuable
comments that allowed us to improve the presentation of our results
and for drawing our attention to reference
\cite{kalashnikov94}. J.-B.~B.\ thanks Agence Nationale de la
Recherche PIECE 12-JS01-0006-01 and Nucleo Milenio NC120062 for
partial support. A.~C.\ thanks support of PUCV Projects 126.711/2014
and 37.375/2014. J.~F.\ was partially supported by Basal-Conicyt,
Nucleo Milenio NC120062 and Fondecyt Project 1150570.

\bibliographystyle{plain}
\bibliography{biblio}

\end{document}